\documentclass[11pt]{article}
\usepackage{amssymb,amsfonts,amsmath, psfrag,eepic,colordvi}
\newtheorem{theo}{Theorem}

\makeatletter \@addtoreset{equation}{section}
\@addtoreset{theo}{}\makeatother

\def\qed{\hfill \rule{4pt}{7pt}}
\def\pf{\noindent {\it Proof.} }

\begin{document}
\begin{center}
{\LARGE\bf On the Hook Length Formula for
Binary Trees}\\[12pt]
William Y.C. Chen$^1$  $\quad$ and $\quad$  Laura L.M. Yang$^2$ \\[6pt]
Center for Combinatorics, LPMC\\
Nankai University,
Tianjin 300071, P. R. China\\[5pt]
$^1$chen@nankai.edu.cn, $^2$yanglm@hotmail.com
\end{center}

\vskip 5mm

\noindent {\bf Abstract. }   We present a simple combinatorial
proof of Postnikov's hook length formula for binary trees.

 \vskip 8mm


 
\vskip 1cm

Let $[n]=\{1,2, \ldots, n\}$.  It is well known that the number of
labeled trees on $[n]$ equals $n^{n-2}$, and the number of rooted
trees on $[n]$ equals $n^{n-1}$ \cite{moon_y1970, stanley2}.
Recently, Postnikov \cite{postnikov} derived the following
identity on binary trees and asked for a combinatorial proof
\cite{postnikov}. We adopt the terminology of Postnikov
\cite{postnikov}. Given a binary tree $T$ and a vertex $v$ of $T$,
we use $h(v)$ to denote the ``hook-length'' of  $v$, namely, the
number of descendants of $v$ (including $v$ itself). Postnikov
\cite{postnikov} obtained the following identity.

\begin{theo} For $n\geq 1$, we have
\begin{equation}\label{binary}
(n+1)^{n-1}=\sum_{T}\frac{n!}{2^n}\prod_{v\in
T}\left(1+\frac{1}{h(v)}\right),
\end{equation}
where the  sum ranges over all binary trees $T$ with $n$ vertices.
\end{theo}

Our combinatorial proof is based on the following equivalent
formulation of (\ref{binary}) in terms of rooted trees:
\begin{equation}\label{binary2}
(n+1)^{n}=\sum_{T}\frac{(n+1)!}{2^n}\prod_{v\in
T}\left(1+\frac{1}{h(v)}\right),
\end{equation}

\pf Let $F_n$ denote the quantity on the right hand side of
(\ref{binary2}). For any unlabeled  binary tree $T$ with  $n$
vertices, the hook length of the root is always $n$. Let us
consider  binary trees $T$ such that the left subtree $T_1$ has
$k$ vertices and the right subtree $T_2$ has $n-k-1$ vertices.
From the relation
\[ {(n+1)!\over 2^n} \left(1+{1\over n}\right) = {n+1 \over 2n} \,{n+1\choose k+1}
 {(k+1)!\over 2^k} \,
{(n-k)!\over 2^{n-k-1}},
\]
we have
\begin{eqnarray*}
F_n &=& \frac{n+1}{2n} \sum_{k=0}^{n-1}{n+1 \choose k+1}
\sum_{T_1}\frac{(k+1)!}{2^k}\prod_{v\in
T_1}\left(1+\frac{1}{h(v)}\right)
\sum_{T_2}\frac{(n-k)!}{2^{n-k-1}} \prod_{v\in
T_2}\left(1+\frac{1}{h(v)}\right),
\end{eqnarray*}
where $T_1$ (or $T_2$) ranges over all binary trees on $k$ (or
$n-k-1$) vertices. Hence $F_n$ satisfies the following recurrence
relation:
\begin{equation}\label{F}
F_n=\frac{n+1}{2n} \sum_{k=0}^{n-1}{n+1 \choose k+1}F_kF_{n-k-1}.
\end{equation}
It is  known that the number $T_n=n^{n-2}$ of labeled trees with
$n$ vertices has  the same recurrence relation:
\begin{equation}\label{tree}
2nT_{n+1}=\sum_{k=0}^{n-1}{n+1 \choose
k+1}(k+1)T_{k+1}(n-k)T_{n-k}.
\end{equation}
Let $R_n=n T_n$ denote the number of rooted tree on $n$ vertices.
Then the above recurrence (\ref{tree}) can be recast as
\begin{equation}\label{T}
R_{n+1}=\frac{n+1}{2n}\sum_{k=0}^{n-1}{n+1 \choose
k+1}R_{k+1}R_{n-k},
\end{equation}
A combinatorial interpretation of (\ref{tree}) is given by
Moon\cite{moon_y1970}: The left hand side of (\ref{tree}) equals
the number of labeled trees on $[n+1]$ with a distinguished edge
and a direction on this distinguished edge. Let $T$ be such a
tree, we may decompose it into an ordered pair of rooted trees by
cutting off the distinguished edge.

Combining the recurrence (\ref{F}) of $F_n$ with the recurrence
(\ref{T}) of $R_n$,  we arrive at the conclusion that $
(n+1)T_{n+1}=F_n$. Hence we obtain (\ref{binary2}). \qed

 S. Seo \cite{seo} also found combinatorial proof of
the identity (\ref{binary}).

\vskip 5mm

\noindent{\bf Acknowledgments.} This work was done under the
auspices of the  ¡°973¡± Project on Mathematical Mechanization,
the National Science Foundation, the Ministry of Education, and
the Ministry of Science and Technology of China.

\small

\end{document}